\def\version{February 23, 2002}

\documentclass[reqno,11pt]{amsart}

\usepackage{amssymb, epsf}

\def\be{\begin{equation}}
\def\ba{\begin{align}}
\def\bm{\begin{multline}}
\def\bfig{\begin{figure}[htb]}
\def\efig{\end{figure}}

\newcommand{\bibit}[1]{\bibitem[#1]{#1}}
\newcommand{\paper}[1]{{\it #1}, }
\newcommand{\journal}[4]{#1 {\bf #2}, #3 (#4)}
\newcommand{\CMP}{Commun.\ Math.\ Phys.}

\newcommand{\PTRF}{Probab.\ Th.\ Rel.\ Fields}
\newcommand{\RMP}{Rev.\ Math.\ Phys.}

\numberwithin{equation}{section}
\newtheorem{theorem}{Theorem}
\newtheorem{proposition}[theorem]{Proposition}
\newtheorem{lemma}[theorem]{Lemma}

\newcounter{eqs}
\newcommand{\startalpheqno}{\setcounter{eqs}{\value{equation}}
\setcounter{equation}{0} \addtocounter{eqs}{1}
\renewcommand{\theequation}{\arabic{section}.\arabic{eqs}{\it \alph{equation}}}
}
\newcommand{\stopalpheqno}{\setcounter{equation}{\value{eqs}}
\renewcommand{\theequation}{\arabic{section}.\arabic{equation}}}

\newcommand{\nn}{\nonumber}

\renewcommand{\leq}{\;\leqslant\;}
\renewcommand{\geq}{\;\geqslant\;}

\newcommand{\dd}{{\rm d}}

\newcommand{\sumtwo}[2]{\sum_{\substack{#1 \\ #2}}}

\newcommand{\prodtwo}[2]{\prod_{\substack{#1 \\ #2}}}

\def\dist{{\operatorname{dist\,}}}

\newcommand{\compl}{{\text{\rm c}}}

\newcommand{\const}{{\text{\rm const}}}
\def\upchi{\raise1pt\hbox{$\chi$}}

\def\writefig#1 #2 #3 {\rlap{\kern #1 truecm \raise #2 truecm
\hbox{#3}}}
\def\figtext#1{\smash{\hbox{#1}} \vspace{-5mm}}

\newcommand{\caC}{{\mathcal C}}

\newcommand{\caG}{{\mathcal G}}

\newcommand{\caL}{{\mathcal L}}

\newcommand{\bbR}{{\mathbb R}}

\newcommand{\bbZ}{{\mathbb Z}}

\newcommand{\bsm}{{\boldsymbol m}}

\begin{document}

{\small
\hfill \version
}
\vspace{2mm}

\title{A self-avoiding walk with attractive interactions}

\author{Daniel Ueltschi}

\address{Daniel Ueltschi,
Department of Mathematics,
University of California,
Davis, CA 95616, USA;
{\rm ueltschi@math.ucdavis.edu.}}

\thanks{Work partially supported by the US National Science
Foundation under grant PHY-98 20650.}

\maketitle

\begin{quote}
{\small
{\bf Abstract.}
A self-avoiding walk with small attractive interactions is described here.
The existence of the connective constant is established, and the diffusive behavior is
proved using the method of the lace expansion.
}

\vspace{1mm}
\noindent
{\footnotesize {\it Keywords:} Self-avoiding random walks, lace expansion.

\noindent
{\it 2000 Math.\ Subj.\ Class.:} 60K35, 60G50, 82B41.}
\end{quote}

\section{Introduction}

A powerful tool for the study of self-avoiding walks is the lace expansion of Brydges and
Spencer \cite{BS}. It is applicable above four dimensions and shows the mean-field behavior
of self-avoiding walks, that is, critical exponents are those of the simple random walk. An
extensive survey of random walks can be found in \cite{MS}.

The lace expansion was originally introduced for weakly self-avoiding walks, and extended
to the fully self-avoiding case by Slade and Hara \cite{Sla,HS}. Several improvements,
simplifications and alternate approaches have since been proposed, see \cite{GI, KLMS, vHHS,
vHS}. A recent work by Bolthausen and Ritzmann \cite{BR} uses a
fixed point argument and avoids the difficulties that are present when working in the
Fourier space. Its actual range of applicability is limited to the case of small repulsions, but
an extension to the self-avoiding case may be possible.

The purpose of this article is to show that the lace expansion can also be used when the walk experiences small
nearest-neighbor attractions. We consider a model of random walks $w = (w_0, \dots, w_n)$ with $w_t \in
\bbZ^d$, where the connectivity $C_n(x)$ between 0 and $x$ is defined by
\be
\label{defCn}
C_n(x) = \sumtwo{w: 0 \to x}{|w|=n} \prod_{t=1}^n D(w_t - w_{t-1}) \prod_{0\leq s<t\leq n}
\bigl( 1 - U(w_s-w_t) \bigr);
\end{equation}
the sum is over all $n$ steps random walks $w \in (\bbZ^d)^{n+1}$ with $w_0=0$ and $w_n=x$.
The jumps of the walk are given a positive weight $D$, that has the symmetries of the
lattice (precisely: invariance under permutations and inversions of coordinates), and that satisfies the following assumptions:
\be
\label{Delta}
D(0)=0, \quad \sum_y D(y)=1, \quad \inf_{|x-y|=1, y\neq0} \frac{D(y)}{D(x)} = \Delta>0.
\end{equation}
Here and in the sequel $|x|$ denotes the $\ell_2$ norm of $x \in \bbZ^d$.
The last ``smoothness" condition is not very desirable, but it plays an important technical role;
notice that it allows a weight with exponential decay. The potential $U$ is
\be
U(x) = \begin{cases} 1 & \text{if } x=0 \\ -\kappa & \text{if } |x|=1 \\ 0 &
\text{otherwise.} \end{cases}
\end{equation}
Here we shall take $\kappa$ small and positive. Let
\be
c_n = \sum_x C_n(x);
\end{equation}
one easily gets $c_1 = 1 + 2d\kappa D(1)$, where $D(1) = D(x)$ with $|x|=1$.
We first establish the existence of the connective constant $\mu = \lim_n c_n^{1/n}$.
Theorem \ref{thmconnconst} is valid for all dimensions and small attractions.

\begin{theorem}
\label{thmconnconst}
Assume that $\kappa$ is small enough, so as to satisfy
$$
(1+\kappa)^{2d} \leq 1 + \frac{\Delta^2}{2d \, (1+\kappa)^{2d-1}}.
$$
Then the sequence $(c_n^{1/n})$ converges to a number $\mu$ with $2^{-d} \leq \mu \leq c_1$.
\end{theorem}

The proof of this theorem is given in Section \ref{secconnconst}.

Next we state a result that will be proved using the lace expansion method. The expansion is
rather easy to perform; the difficult task is to prove the convergence. This will be done in two
steps. First, we shall obtain bounds on lace expansion terms involving the supremum norm of
$C_n(x)$. Second, we shall check the hypothesis of van der Hofstad and Slade \cite{vHS};
their results imply Theorem \ref{thmdiffbeh} below.

We consider a positive differentiable even function $h(\xi)$ on $\bbR$, that is decreasing for
$\xi>0$, and that satisfies
\startalpheqno
\ba
\label{condha}
&\int |\xi|^{d+1+3\varepsilon} h(\xi) \dd\xi < \infty \quad \text{for some } \varepsilon \in (0, 1
\wedge \tfrac{d-4}4), \\
&\sup_{\xi \in \bbR} \; \Bigl| \frac{h'(\xi)}{h(\xi)} \Bigr| < \infty.
\label{condhb}
\end{align}
\stopalpheqno
For $x\in\bbZ^d \setminus \{0\}$, we define
\be
\label{defD}
D(x) = \frac{h(|x|/L)}{\sum_{y\in\bbZ^d \setminus \{0\}} h(|y|/L)},
\end{equation}
and $D(0)=0$. The condition \eqref{condha} is a technical one that appears in
\cite{vHS}; \eqref{condhb} ensures the existence of a non-zero constant $\Delta$, see
\eqref{Delta}, at least when $L$ is large.

\begin{theorem}
\label{thmdiffbeh}
Suppose $d\geq5$ and define $D$ by \eqref{defD} with $h$ satisfying \eqref{condha} and
\eqref{condhb}. There exists $L_0<\infty$ such that if $L \geq L_0$, and if $\kappa$ is small
enough so that the condition of Theorem \ref{thmconnconst} holds true, the mean-square displacement satisfies
$$
\frac1{c_n} \sum_x |x|^2 C_n(x) = n \delta \bigl[ 1 + O(n^{-\varepsilon}) \bigr].
$$
The diffusion constant $\delta$ can be given an explicit expression, see \eqref{diffconst}.
\end{theorem}

A self-avoiding random walk with strong attractions ($\kappa$ large) displays a very
different behavior. A typical walk is expected to maximize nearest-neighbor contacts and
to occupy as little a space as possible. We can actually compute a lower bound for
the connective constant by considering only such walks. Let $\caC(n^{1/d})$ denote the cube
of size $n^{1/d}$ centered at the origin, and define
\be
\gamma = \liminf_{n\to\infty} \Bigl[ \sumtwo{w \subset \caC(n^{1/d})}{|w|=n} \prod_{t=1}^n D(w_t
- w_{t-1}) \Bigr]^{1/d}.
\end{equation}
The sum is over all self-avoiding walks starting at the origin and with support
$\caC(n^{1/d})$. We easily obtain
\be
\mu \geq \gamma (1+\kappa)^d.
\end{equation}
Since $\gamma>D(1)$, we see that the bound $\mu \leq c_1$ given in Theorem
\ref{thmconnconst} cannot be true for $\kappa$ large. One should also expect that the
mean-square displacement has leading term $n^{2/d}$, that is, the critical exponent is
smaller than in the case of small $\kappa$.

A random walk with both on-site repulsion and nearest-neighbor attraction is studied in
\cite{vHK}.

The lace expansion is explained in Section \ref{seclexp}, and suitable bounds
of lace expansion terms are obtained. The special difficulties associated with attractive interactions
are treated with the help of Section \ref{secconnconst}. Section \ref{secdiffbeh} contains the proof of
Theorem \ref{thmdiffbeh}, based on Section \ref{seclexp} and \cite{vHS}. Notice that the
assumptions of \cite{vHS} are less restrictive and the claims are stronger. More
general walks can be considered and a local central limit theorem holds true. See \cite{vHS}
for more informations.

\section{The connective constant}
\label{secconnconst}

This section is devoted to the proof of Theorem \ref{thmconnconst}.
A lower bound for $c_n^{1/n}$ can be found by restricting the sum in \eqref{defCn} to random walks which jump only in
positive directions, and by neglecting the nearest-neighbor attractions. We get
\be
c_n \geq \Bigl( \sum_{x: x_i \geq 0} D(x) \Bigr)^n \geq 2^{-dn};
\end{equation}
we used the fact that $D(x)$ is normalized, and its sum in the first octant is at least
$2^{-d}$. 

We prove below that
\be
\label{subadd}
c_{m+n} \leq c_m c_n.
\end{equation}
Then $c_n \leq c_1^n$, and we obtain the upper bound. From \eqref{subadd} and a standard
subadditivity argument, we get the convergence of the sequence $(c_n^{1/n})$.

The difficulty is to prove \eqref{subadd}. It clearly holds in the case of repulsive
interactions, and fails when only attractions are present. Here, one has to play the
attractions against the self-avoidance, to see that the effective behavior is indeed
repulsive.

Let us introduce
\be
W(w) = \prod_{t=1}^{|w|} D(w_t-w_{t-1}) \prod_{0\leq s<t\leq |w|} \bigl( 1 - U(w_s-w_t)
\bigr);
\end{equation}
then
\be
\label{cmn}
c_{m+n} = \sum_{x,y} \sumtwo{w: 0\to x}{|w|=m} W(w) \sumtwo{w': x\to y}{|w'|=n} W(w')
\prodtwo{0\leq s<m}{0<t\leq n} \bigl( 1 - U(w_s - w_t') \bigr).
\end{equation}

Let us fix $w$; we show that the following holds true for all $0\leq j<m$:
\be
\label{keyineq}
\sumtwo{w': x\to y}{|w'|=n} W(w') \prodtwo{j\leq s<m}{0<t\leq n} \bigl( 1 -
U(w_s-w_t') \bigr) \leq \sumtwo{w': x\to y}{|w'|=n} W(w') \prodtwo{j+1\leq s<m}{0<t\leq n} \bigl( 1 -
U(w_s-w_t') \bigr).
\end{equation}
Notice that $(w_0, \dots, w_{j-1})$ does not play any role in the expression above.
Inequality \eqref{keyineq} allows to remove the product term in \eqref{cmn}, and one
obtains \eqref{subadd}.

\bfig
\epsfxsize=50mm
\centerline{\epsffile{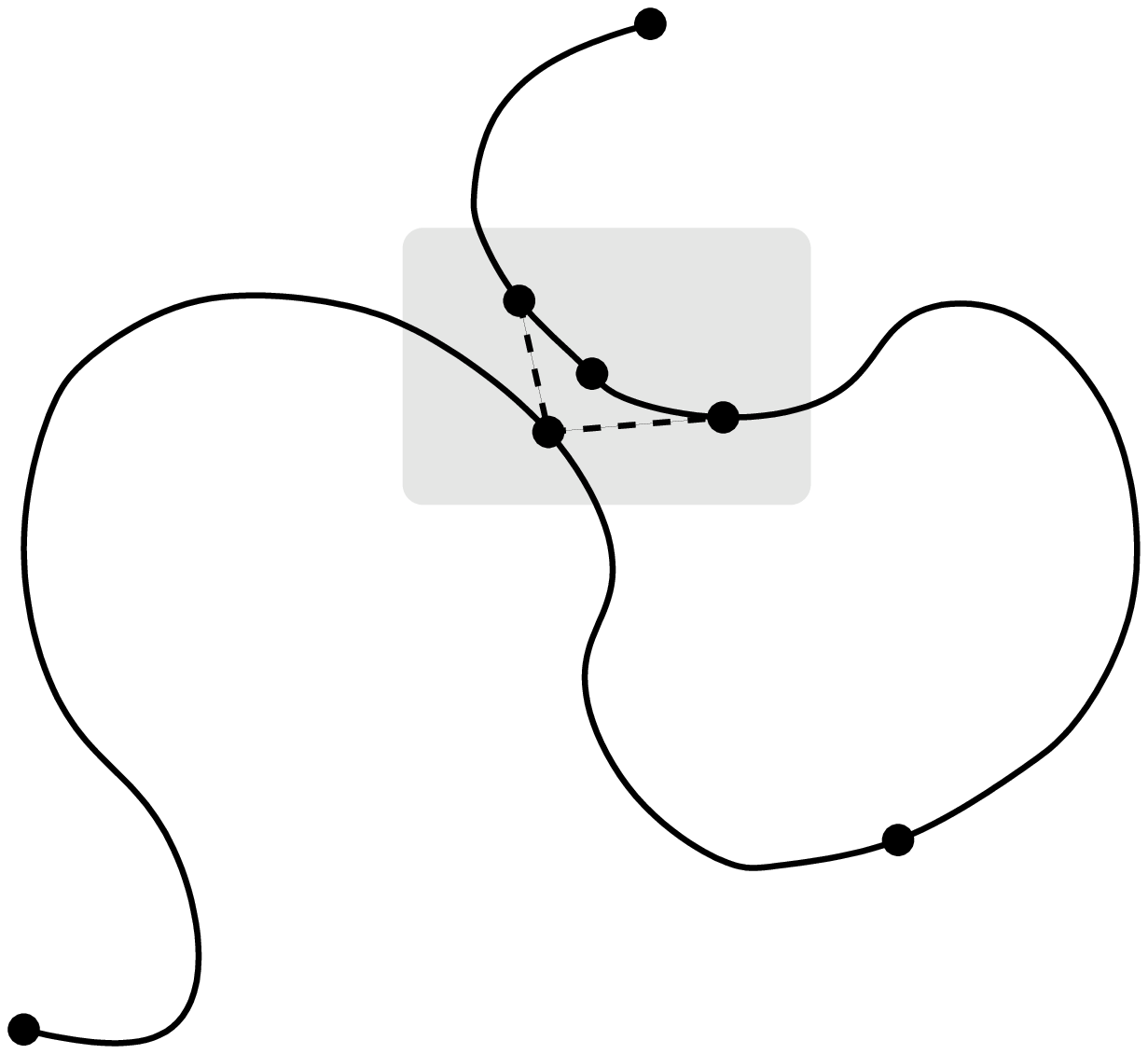}}
\figtext{
\writefig	-2.65	0.6	{\footnotesize 0}
\writefig	1.45	1.2	{\footnotesize $x$}
\writefig	0.5	4.8	{\footnotesize $y$}
\writefig	-2.7	2.7	{$w$}
\writefig	2.45	3.0	{$w'$}
\writefig	-0.55	3.0	{\footnotesize $w_j$}
\writefig	0.2	3.45	{\footnotesize $w_u'$}
}
\caption{The walk $\check w$ defined in \eqref{defwhacek} is a
little deformation of $w'$.}
\label{figdeform}
\end{figure}

Let $\Omega_{xy}^n$ be the set of $n$-steps walks from $x$ to $y$; furthermore, for given
$w_j$ we set
\ba
\Omega_0 &= \{ w' \in \Omega_{xy}^n : \dist(w',w_j) = 0 \}, \nn\\
\Omega_1 &= \{ w' \in \Omega_{xy}^n : \dist(w',w_j) = 1 \}, \nn\\
\Omega_2 &= \{ w' \in \Omega_{xy}^n : \dist(w',w_j) > 1 \}, \nn
\end{align}
where the distance between a walk $w'$ and a point $z$ is $\dist(w',z) = \min_{1\leq t\leq
|w'|} |w_t'-z|$. Clearly,
\be
\sum_{w'\in\Omega_2} W(w') \prodtwo{j\leq s<m}{0<t\leq n} \bigl( 1 -
U(w_s-w_t') \bigr) = \sum_{w'\in\Omega_2} W(w') \prodtwo{j+1\leq s<m}{0<t\leq n} \bigl( 1 -
U(w_s-w_t') \bigr).
\end{equation}
We turn now to the walks of $\Omega_0 \cup \Omega_1$. For $w' \in \Omega_1$, we define
\be
u = \min\{t\geq1: \dist(w_t',w_j) = 1\}.
\end{equation}
Then we consider a walk $\check w \in \Omega_0$ that is a little deformation of $w'$, namely
\be
\label{defwhacek}
\check w_t = \begin{cases} w_t' & \text{if } t \neq u, \\ w_j & \text{if } t=u; \end{cases}
\end{equation}
this is illustrated in Fig.\ \ref{figdeform}. Notice that
\be
\prod_{t=1}^n D(\check w_t - \check w_{t-1}) \geq \Delta^2 \prod_{t=1}^n D(w_t' -
w_{t-1}').
\end{equation}
The walk $\check w$ may involve less nearest neighbor contacts with itself or $w$, than the
walk $w'$. The difference is no more than $2d-1$ (a consequence of self-avoidance), so that
\bm
\prod_{0\leq s<t\leq n} \bigl( 1 - U(\check w_s - \check w_t) \bigr) \prodtwo{j+1\leq
s<m}{0<t\leq n} \bigl( 1 - U(w_s - \check w_t) \bigr) \\
\geq (1+\kappa)^{-(2d-1)} \prod_{0\leq
s<t\leq n} \bigl( 1 - U(w_s'-w_t') \bigr) \prodtwo{j+1\leq s<m}{0<t\leq n} \bigl( 1 -
U(w_s-w_t') \bigr).
\label{uneborne}
\end{multline}

The right side of \eqref{keyineq} involves walks both of
$\Omega_0$ and $\Omega_1$, while the left side involves only walks of $\Omega_1$.
To each walk $w' \in \Omega_1$ corresponds a walk $\check w \in \Omega_0$,
and the weight of $\check w$ is
bounded below by the weight of $w'$, up to a factor $\Delta^2/(1+\kappa)^{2d-1}$.
No more than $2d$ walks $w' \in \Omega_1$ are mapped on a same $\check w$. Starting with
the right side of \eqref{keyineq}, we can write
\ba
\sum_{w' \in \Omega_1} &W(w') \prodtwo{j+1\leq s<m}{0<t\leq
n} \bigl( 1 - U(w_s-w_t') \bigr) + \sum_{\check w \in \Omega_0} W(\check w) \prodtwo{j+1
\leq s<m}{0<t\leq
n} \bigl( 1 - U(w_s - \check w_t) \bigr) \nn\\
&\geq \sum_{w' \in \Omega_1} W(w') \prodtwo{j+1\leq s<m}{0<t\leq
n} \bigl( 1 - U(w_s-w_t') \bigr) \nn\\
&\hspace{3cm} + \frac{\Delta^2}{2d \, (1+\kappa)^{2d-1}} \sum_{w' \in
\Omega_1} W(w') \prodtwo{j+1\leq s<m}{0<t\leq
n} \bigl( 1 - U(w_s-w_t') \bigr) \nn\\
&\geq \frac1{(1+\kappa)^{2d}} \Bigl[ 1 + \frac{\Delta^2}{2d \, (1+\kappa)^{2d-1}} \Bigr]
\sum_{w' \in \Omega_1} W(w') \prodtwo{j\leq s<m}{0<t\leq n} \bigl( 1 -
U(w_s-w_t') \bigr). \nn
\end{align}
The assumption of the theorem implies that the factor in front of the last sum is larger than 1.

The importance of the ``smoothness" condition for $D$ is clear from the occurrence of the
constant $\Delta$ in the equation above. Self-avoidance allowed to write the
inequality \eqref{uneborne}. In the case of weakly self-avoiding walks some sites receive
many visits, and the method described here does not work --- Eq.\ \eqref{keyineq}
actually ceases to be true. While weakly self-avoiding walks should also display
effective repulsion, to prove it looks difficult.

\section{The lace expansion}
\label{seclexp}

The goal now is to write down a lace expansion for our
self-avoiding walk with attractive nearest-neighbor interactions, and then to prove a key
estimate; see Proposition \ref{propestim} below. It will be used in showing the convergence of the
lace expansion, hence in establishing the diffusive behavior of the walk.

The first step consists in obtaining an expansion for
the connectivity $C_n(x)$.
A natural idea is to proceed as in a cluster expansion and, in \eqref{defCn}, to expand the product over
$(s,t)$ so as to get a sum over graphs of $n+1$ vertices, and then to attempt to control
the resulting terms. Dealing with these terms is no easy task, but Brydges and Spencer
have shown that suitable bounds can indeed be proven \cite{BS}. The idea is to take advantage of the
one-dimensional nature of a walk. We consider graphs whose sets of vertices are intervals $[a,b]$ in $\bbZ$; we write $\caG[a,b]$ for the set of all graphs on $[a,b]$, and $\caC[a,b]$ for the set of {\it connected} graphs: a graph $G$ is connected iff
\begin{itemize}
\item both $a$ and $b$ are endpoints of edges of $G$;
\item $\forall c \in (a,b)$: $\exists st \in G$ such that $s<c<t$.
\end{itemize}

A {\it lace} is a minimally connected graph, i.e.\ a connected graph such that the removal of
any edge results in a disconnected graph. We denote by $\caL[a,b]$ the set of laces on
$[a,b]$. If $G \in \caC[a,b]$, one can obtain a lace $L(G) \subset G$ by keeping edges
$s_1t_1, \dots, s_mt_m$ of $G$, according to the following rule:
\begin{itemize}
\item $s_1=a$, $t_1 = \max\{ t: (a,t) \in G \}$
\item $t_2 = \max\{ t: \exists s<t_1 \text{ such that } st \in G \}$, $s_2 = \min\{ s: st_2 \in G \}$
\item[] \hspace{3cm} $\vdots$
\item $t_m=b$, $s_m = \min\{ s: sb \in G \}$.
\end{itemize}

\bfig
\epsfxsize=120mm
\centerline{\epsffile{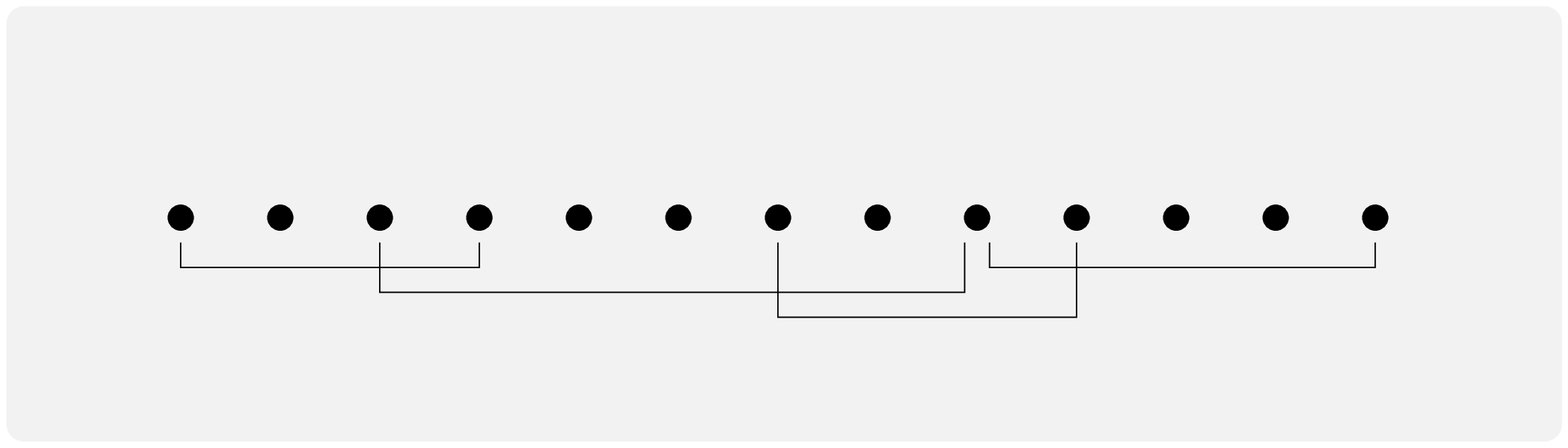}}
\figtext{
\writefig	-4.55	2.55	{\footnotesize $a$}
\writefig	4.35	2.55	{\footnotesize $b$}
\writefig	-4.6	1.7	{\footnotesize $s_1$}
\writefig	-3.1	1.55	{\footnotesize $s_2$}
\writefig	-0.15	1.35	{\footnotesize $s_3$}
\writefig	-2.35	1.5	{\footnotesize $t_1$}
\writefig	1.0	2.55	{\tiny $t_2=s_4$}
\writefig	2.15	1.3	{\footnotesize $t_3$}
\writefig	4.35	1.7	{\footnotesize $t_4$}
}
\caption{Example of a lace. Note that this graph would not be connected without the
presence of $s_3t_3$.}
\label{figlace}
\end{figure}

Let $L$ be a lace. An edge $st \notin L$ such that the lace corresponding to $L \cup \{st\}$
is $L$, is said to be {\it compatible} with $L$, and we write $st \sim L$. Any graph $G$
such that $L(G) = L$ contains all edges of $L$, and edges that are compatible with $L$ (and
reciprocally).

We are looking for an induction relation for $C_n(x)$. We start by rewriting \eqref{defCn}
as
\be
C_n(x) = \sumtwo{w: 0 \to x}{|w|=n} D(w) \sum_{G \in \caG[0,n]} \prod_{st
\in G} \bigl( -U(w_s-w_t) \bigr),
\end{equation}
with
\be
D(w) = \prod_{t=1}^{|w|} D(w_t-w_{t-1}).
\end{equation}

Some graphs have edges attached to 0, and some graphs do not have such edges. Graphs of the
former class can be split into a connected graph containing 0, and another graph whose
support consists of the remaining sites. This leads to the decomposition
\bm
\sum_{G \in \caG[0,n]} \prod_{st \in G} \bigl( -U(w_s-w_t) \bigr) = \sum_{G
\in \caG[1,n]} \prod_{st \in G} \bigl( -U(w_s-w_t) \bigr) \\
+ \sum_{m=1}^n \sum_{G \in \caC[0,m]} \prod_{st \in G} \bigl( -U(w_s-w_t) \bigr)
\sum_{G' \in \caG[m+1,n]} \prod_{st \in G'} \bigl( -U(w_s-w_t) \bigr).
\end{multline}
Then
\bm
C_n(x) = D * C_{n-1}(x)
+ \sumtwo{w: 0 \to x}{|w|=n} D(w) \sum_{m=1}^n \sum_{L \in \caL[0,m]}
\prod_{st \in L} \bigl( -U(w_s-w_t) \bigr) \\
\prod_{st \sim L} \bigl( 1 -
U(w_s-w_t) \bigr) \prod_{m < s < t \leq n} \bigl( 1 - U(w_s-w_t)
\bigr).
\end{multline}
The star symbol denotes the convolution of $D$ and $C_{n-1}$, namely $\sum_y D(y) C_{n-1}(x-y)$.
We define
\be
\Pi_m(x) = \sumtwo{w: 0 \to x}{|w|=m} D(w) \sum_{L \in \caL[0,m]} \prod_{st \in L}
\bigl( -U(w_s-w_t) \bigr) \prod_{st \sim L} \bigl( 1 - U(w_s-w_t) \bigr),
\end{equation}
and $\pi_m = \sum_x \Pi_m(x)$. Notice that $\pi_1 = 2d\kappa D(1)$. Setting $C_0(x) =
\delta_{0x}$, we get the desired formula:
\be
C_n(x) = D * C_{n-1}(x) + \sum_{m=1}^n \Pi_m * C_{n-m}(x).
\label{recCn}
\end{equation}

Such a relation is true for simple random walks, setting $\Pi_m(y) \equiv 0$. The second term
is therefore the correction due to the self-interactions, and the purpose of the lace expansion
is to show that it is small.

Let $\caL^{(N)}[a,b]$ denote the set of laces on $[a,b]$ with exactly $N$ edges. We write
\be
\Pi_n(x) = \sum_{N \geq 1} \Pi_n^{(N)}(x),
\end{equation}
with
\be
\label{defPinN}
\Pi_n^{(N)}(x) = \sumtwo{w: 0\to x}{|w|=n} D(w) \sum_{L \in \caL^{(N)}[0,n]} \prod_{st \in L}
\bigl( -U(w_s-w_t) \bigr) \prod_{st \sim L} \bigl( 1 - U(w_s-w_t) \bigr).
\end{equation}

In order to prove the convergence of the lace expansion, one needs bounds on $\Pi_m(x)$. We
propose here estimates that involve norms of $C_n$, and norms of moments of $C_n$. They are both standard and
useful. The following proposition holds true in all dimensions, and
with all $D$ satisfying \eqref{Delta}. It will only be used in the restricted situation of
Theorem \ref{thmdiffbeh}, however.

\begin{proposition}
\label{propestim}
If $\kappa$ is small enough so as to satisfy the condition in Theorem \ref{thmconnconst},
we have the following bounds:
\begin{itemize}
\item[(i)] For $N=1$,
$$
\|\Pi_n^{(1)}\|_1 \leq (1+2d\kappa) \|C_{n-1}\|_\infty.
$$
\item[(ii)] For $N\geq2$,
$$
\|\Pi_n^{(N)}\|_1 \leq (2N-1) 2^{N-1} (1+2d\kappa)^N \sum \prodtwo{j=1}{\rm odd}^{2N-1}
\|C_{m_j}\|_\infty \prodtwo{j=2}{\rm even}^{2N-2} \|C_{m_j}\|_1.
$$
\item[(iii)] For $N=1$ and all $\gamma>0$,
$$
\bigl\| |x|^\gamma \Pi_n^{(1)} \bigr\|_1 \leq 2d\kappa \|C_{n-1}\|_\infty.
$$
\item[(iv)] For $N\geq2$ and all $1\leq\gamma\leq2$,
\ba
\bigl\| |x|^{2\gamma} \Pi_n^{(N)} \bigr\|_1 \leq &(N-1)^{2\gamma-2} (2N-1) 2^{2\gamma-2+N} 
(1+2d\kappa)^N \nn\\
&\sum \sumtwo{i=3}{\rm odd}^{2N-1} \bigl\| [|x|^\gamma+1] C_{m_i} \bigr\|_\infty
\prodtwo{i'=1, i'\neq i}{\rm odd}^{2N-1} \|C_{m_{i'}}\|_\infty \nn\\
&\sumtwo{j=2}{\rm even}^{2N-2} \bigl\| [|x|^\gamma+1] C_{m_j} \bigr\|_1 \prodtwo{j'=2,
j'\neq j}{\rm even}^{2N-2} \|C_{m_{j'}}\|_1. \nn
\end{align}
\end{itemize}
Unlabeled sums in {\rm (ii)} and {\rm (iv)} are over $m_1, \dots, m_{2N-1}$ whose sum is $n$, and such that $m_1$ is the
larger number, and $m_{2j} \leq m_{2j+1}$ for all $1\leq j\leq N-1$.
\end{proposition}

\begin{proof}
The proof is standard, except for the difficulties associated with the attractive
interactions. For part (i),
\ba
\|\Pi_n^{(1)}\|_1 &\leq \sum_x \sumtwo{w:0\to x}{|w|=n} |U(x)| D(w) \prodtwo{0\leq s<t\leq
n}{(s,t)\neq(0,n)} \bigl( 1-U(w_s-w_t) \bigr) \nn\\
&= \sum_x |U(x)| \sum_y D(y) \sumtwo{w:y\to x}{|w|=n-1} W(w) \prod_{t=0}^{n-1} \bigl(
1-U(w_t) \bigr).
\label{bornePi1}
\end{align}
This is is a special case of \eqref{keyineq}: in \eqref{keyineq}, take $m=1$, $j = 0$, and
$w_j = 0$. As a result, we get an upper bound by removing the product in
\eqref{bornePi1}, and we easily obtain Proposition \ref{propestim} (i).

\bfig
\epsfxsize=120mm
\centerline{\epsffile{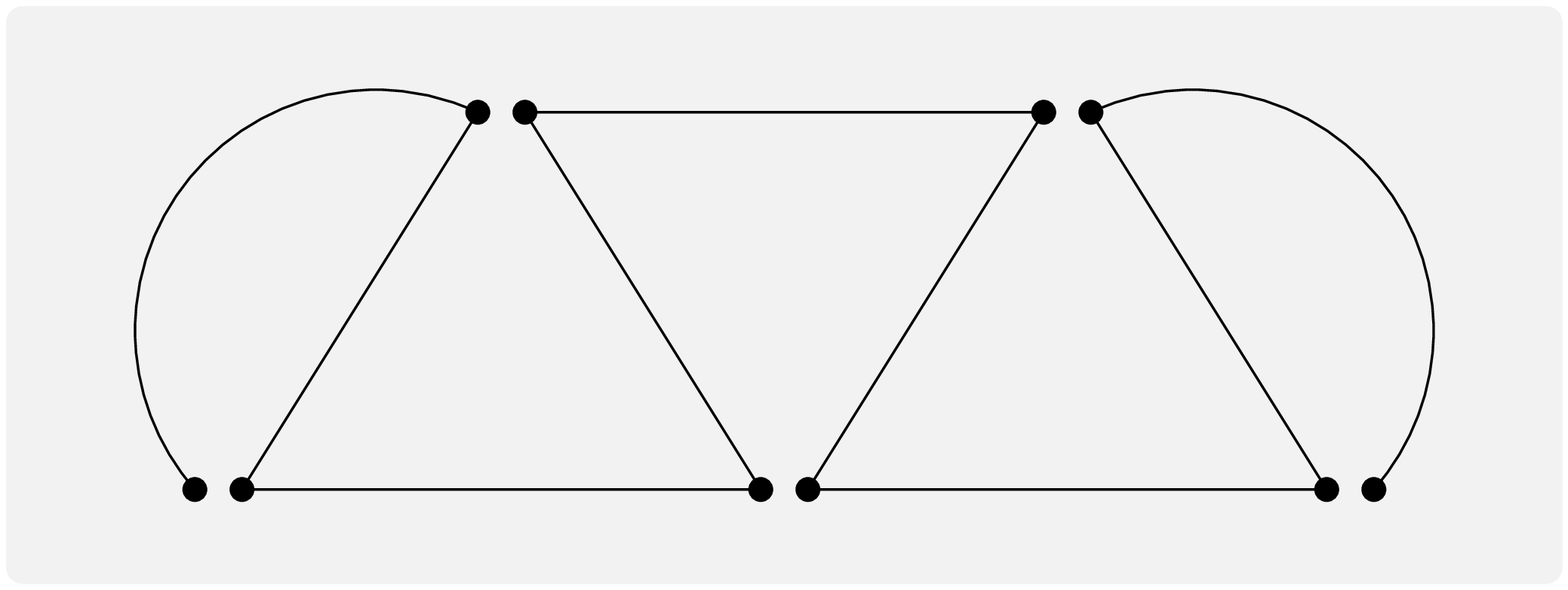}}
\figtext{
\writefig	-4.5	0.98	{\footnotesize $0$}
\writefig	-4.15	1.0	{\footnotesize $x_0'$}
\writefig	-2.45	4.35	{\footnotesize $x_1$}
\writefig	-1.95	4.35	{\footnotesize $x_1'$}
\writefig	-0.4	1.0	{\footnotesize $x_2$}
\writefig	0.05	1.0	{\footnotesize $x_2'$}
\writefig	1.75	4.35	{\footnotesize $x_3$}
\writefig	2.15	4.35	{\footnotesize $x_3'$}
\writefig	3.8	1.0	{\footnotesize $x_4$}
\writefig	4.25	1.0	{\footnotesize $x_4'$}
\writefig	-5.0	3.1	{\footnotesize $1$}
\writefig	-3.3	2.9	{\footnotesize $2$}
\writefig	-2.25	1.45	{\footnotesize $3$}
\writefig	-1.25	2.5	{\footnotesize $4$}
\writefig	-0.15	3.9	{\footnotesize $5$}
\writefig	0.8	2.8	{\footnotesize $6$}
\writefig	1.95	1.45	{\footnotesize $7$}
\writefig	2.8	2.6	{\footnotesize $8$}
\writefig	4.35	3.1	{\footnotesize $9$}
}
\caption{Diagram for $\Pi_n^{(N)}$ with $N=5$. $0, x_1, x_0', x_2, \dots, x_4'$ are successive positions of the
walks (and are summed upon). The legs 1,\dots,9 have length $m_1, \dots, m_9$
respectively, with $m_1 + \dots + m_9 = n$.}
\label{figdiag}
\end{figure}

Figure \ref{figdiag} depicts a diagram that represents the sum over laces in \eqref{defPinN}. A lace is completely determined by $\bsm =
(m_1, \dots, m_{2N-1})$ such that $m_1 + \dots + m_{2N-1} = n$. These intervals satisfy moreover $m_1 \geq
1$, $m_{2N-1} \geq 1$, and for $1 \leq j \leq N-1$: $m_{2j} \geq 1$, $m_{2j+1} \geq 0$.
Then (with $x_0 = 0$)
\bm
\|\Pi_n^{(N)}\|_1 \leq \sum_\bsm \sumtwo{x_1, \dots, x_{N-1}}{x_0', \dots, x_{N-1}'}
\prod_{j=0}^{N-1} |U(x_j-x_j')| \sumtwo{w^{(1)}: 0 \to
x_1}{|w^{(1)}| = m_1} \; \sumtwo{w^{(2)}: x_1 \to x_0'}{|w^{(2)}| = m_2} \dots \\
\dots \sumtwo{w^{(2N-2)}: x_{N-1}
\to x_{N-2}'}{|w^{(2N-2)}| = m_{2N-2}} \; \sumtwo{w^{(2N-1)}: x_{N-2}'
\to x_{N-1}'}{|w^{(2N-1)}| = m_{2N-1}} \; \prod_{j=1}^{2N-1} D(w^{(j)})
\prod_{st \sim L} \bigl( 1 - U(w_s-w_t) \bigr)
\label{bound}
\end{multline}
where the last product is over all edges compatible with the lace $L$, that is defined by
$\bsm$. The walk $w$ that appears in the last product is the union (`concatenation') of the
walks $w^{(1)}, \dots, w^{(2N-1)}$. All edges
between vertices of a same leg are compatible with $L$, and therefore appear in the product.

We need to get rid of the interactions between different legs. When the random walk is
only repulsive this is easy: neglecting these interactions yields an upper bound. Here we
proceed as in the proof of the existence of the connective constant, using the fact that
the legs are effectively repulsive.

We start with the first leg. It interacts with the legs 2,3,4 only. Notice that all edges
with one endpoint on the leg 1, and the other endpoint on leg 2, 3, or 4, are compatible with
$L$. Let $w$ be a walk on the time interval $[m_1,n]$. Setting $m' =
m_1+m_2+m_3+m_4$, we have that for all $j\geq1$,
\bm
\sumtwo{w^{(1)}: 0\to x_1}{|w^{(1)}|=m_1} W(w^{(1)}) \prodtwo{st \sim
L: 0 \leq s<m_1}{m_1+j \leq t < m'} \bigl( 1 - U(w_s^{(1)} - w_t)
\bigr) \\
\leq \sumtwo{w^{(1)}: 0\to x_1}{|w^{(1)}|=m_1} W(w^{(1)})
\prodtwo{st \sim L: 0\leq s<m_1}{m_1+j+1 \leq t < m'} \bigl( 1 - U(w_s^{(1)} - w_t) \bigr).
\label{bonneborne}
\end{multline}
Indeed, the restriction of $w$ to the time interval $[m_1,m']$ is a self-avoiding walk, and
we are therefore in the same situation as \eqref{keyineq}.

Inequality \eqref{bonneborne} implies that we get an upper bound by neglecting the
interactions between the first leg and the others in \eqref{bound}. The second leg interacts
only with legs 3 and 4, and a similar inequality can be written. We proceed further by
considering the interactions between the third leg and the subsequent ones (precisely: the
legs 4,5,6), and so on. At the end we have an upper bound by removing all interactions
between different legs. Hence,
\bm
\|\Pi^{(N)}_n\|_1 \leq \sum_\bsm \sumtwo{x_1, \dots, x_{N-1}}{x_0', \dots, x_{N-1}'}
\prod_{j=0}^{N-1} |U(x_j-x_j')| \,
C_{m_1}(x_1) C_{m_2}(x_1-x_0') \\
\Bigl( \prod_{j=2}^{N-1} C_{m_{2j-1}}(x_j - x_{j-2}') C_{m_{2j}}(x_j - x_{j-1}') \Bigr)
C_{m_{2N-1}}(x_{N-1}' - x_{N-2}').
\label{bound2}
\end{multline}

The rest of the proof is standard.
We first sum over the edge with the larger $m$. Let us denote the corresponding index by $\ell$. Then
we group remaining edges into pairs, see Figure \ref{figpairing} for examples. We define $J=J(\bsm)$ to
be the set of indices such that $\ell \in J$, and if $(i,j)$ denotes paired edges, then $i \in J$ if
$m_i \geq m_j$, and $j \in J$ otherwise. We
obtain
\be
\|\Pi_n^{(N)}\|_1 \leq \sum_\bsm \prod_{j \in J} \|C_{m_j}\|_\infty
\sumtwo{x_1,
\dots, x_{N-1}}{x_0',\dots,x_{N-1}'} \prod_{j=0}^{N-1} |U(x_j-x_j')| \prod_{j \in J^\compl}
C_{m_j}(y_j-y_j').
\label{bound3}
\end{equation}
Here, $y_j$ and $y_j'$ are the endpoints of the leg $j$; they are determined unambiguously by $\bsm$ and
by $x_1, \dots, x_{N-1}'$.
\bfig
\epsfxsize=150mm
\centerline{\epsffile{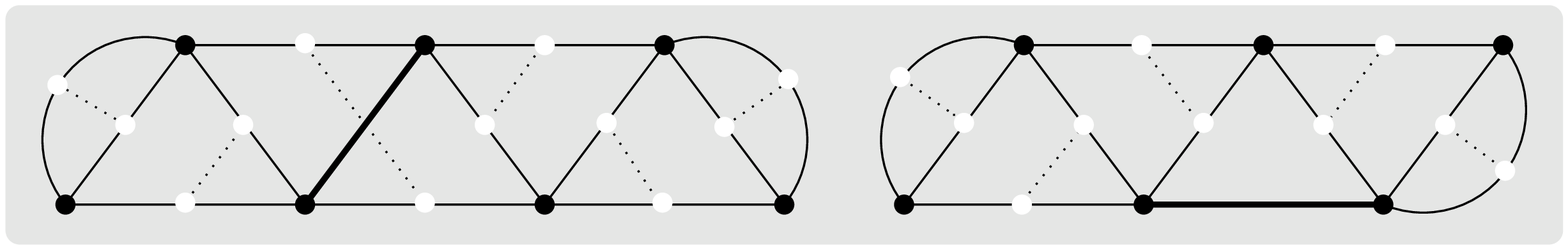}}
\figtext{
\writefig	-4.35	1.6	{\footnotesize $\ell$}
\writefig	4.35	1.15	{\footnotesize $\ell$}
}
\caption{Illustration for the pairing of edges.}
\label{figpairing}
\end{figure}
The pairing of edges was made in such a way that the graph with $N$ vertices, and edges given by
$J^\compl$, is always connected (it is actually a tree). Therefore
\be
\sumtwo{x_1, \dots, x_{N-1}}{x_0',\dots,x_{N-1}'} \prod_{j=0}^{N-1} |U(x_j-x_j')| \prod_{j \in
J^\compl} C_{m_j}(y_j-y_j') \leq (1+2d\kappa)^N \prod_{j \in
J^\compl} \|C_{m_j}\|_1.
\end{equation}
Since there are $(2N-1)$
possibilities for $\ell$, and $2^{N-1}$ for $J \setminus \{\ell\}$, we obtain the bound of Proposition
\ref{propestim} (ii); indeed, the latter corresponds to the case $\ell=1$ and $J = \{ 2j-1
: 1\leq j\leq N \}$.

The bound (iii) is similar to (i), since $\Pi_n^{(N)}(x)=0$ if $|x|>1$.

The proof of (iv) can be done by first
modifying equations \eqref{bound}, \eqref{bound2}, and \eqref{bound3}, turning the left side into
$\sum_x |x|^{2\gamma} |\Pi_n^{(N)}(x)|$, and inserting an extra factor $|x_{N-1}'|^{2\gamma}$  next to the sum over
sites $x_1, \dots, x_{N-1}'$. The self-interaction has range 1, so that $|x_j-x_j'|\leq1$;
introducing appropriate vectors $e_j$ with $|e_j|=0$ or 1, we have
\be
|x_{N-1}'|^\gamma = \Bigl| \sum_{j \in J^\compl} (y_j - y_j' - e_j) \Bigr|^\gamma.
\end{equation}
By H\"older and since $|J^\compl|=N-1$, we get
\ba
|x_{N-1}'|^\gamma &\leq (N-1)^{\gamma-1} \sum_{j \in J^\compl} |y_j-y_j'-e_j|^\gamma \nn\\
&\leq 2^{\gamma-1} (N-1)^{\gamma-1} \sum_{j \in J^\compl} \Bigl[ |y_j-y_j'|^\gamma + 1
\Bigr].
\end{align}
This inequality also holds when $J^\compl$ is replaced by $J\setminus\{\ell\}$. It turns out
that a suitable bound for $|x_{N-1}'|^{2\gamma}$ is the product of bounds with $J^\compl$ and
$J\setminus\{\ell\}$. We obtain
\ba
\sum_x |x|^{2\gamma} |\Pi_n^{(N)}(x)| &\leq 2^{2\gamma-2} (N-1)^{2\gamma-2} \sum_\bsm
\sumtwo{x_1, \dots, x_{N-1}}{x_0',\dots,x_{N-1}'} \prod_{j=0}^{N-1} |U(x_j-x_j')| \nn\\
&\sum_{i \in J\setminus\{\ell\}} \bigl[ |y_i-y_i'|^\gamma + 1 \bigr]
\sum_{j \in J^\compl} \bigl[ |y_j-y_j'|^\gamma + 1 \bigr] C_{m_1}(x_1) C_{m_2}(x_1-x_0') \\
&\Bigl( \prod_{j'=2}^{N-1} C_{m_{2j'-1}}(x_{j'} - x_{j'-2}') C_{m_{2j'}}(x_{j'} - x_{j'-1}') \Bigr)
C_{m_{2N-1}}(x_{N-1}' - x_{N-2}'). \nn
\end{align}
The rest of the proof of item (iv) is similar to (ii).
\end{proof}

\section{The diffusive behavior}
\label{secdiffbeh}

Convergence of the lace expansion follows from Proposition \ref{propestim}, but it is still
a difficult and intricate task. A rather general context was considered in \cite{vHS} that
applies here. The starting point is the following equation,
\be
f_n(k;z) = \sum_{m=1}^n g_m(k;z) f_{n-m}(k;z).
\label{geneq}
\end{equation}
Here, $f_0(k;z)=1$ and $k \in [-\pi,\pi]^d$; $z$ is a positive parameter. An extra term is
allowed in \cite{vHS}, but it is not needed here. One comes close to this equation by
taking the Fourier transform of \eqref{recCn}, namely,
\be
\hat C_n(k) = \hat D(k) \hat C_{n-1}(k) + \sum_{m=1}^n \hat\Pi_m(k) \hat C_{n-m}(k).
\label{Ftransf}
\end{equation}
Comparing with \eqref{geneq}, we see that the $m=1$ term does not perfectly match. Notice that
$\Pi_1(x) = \kappa D(1) \delta_{|x|,1}$ and $C_1(x) = D(x) + \Pi_1(x)$, where $D(1)=D(x)$ with $|x|=1$. One obtains
\eqref{geneq} with the following definitions:
\startalpheqno
\ba
&\hat E(k) = \frac{\hat D(k) + 2\kappa D(1) \sum_{i=1}^d \cos k_i}{1 + 2d\kappa D(1)}, \\
&f_0(k;z) = 1, \quad f_1(k;z) = g_1(k;z) = z \hat E(k), \\
&f_n(k;z) = \Bigl( \frac z{1+2d\kappa D(1)} \Bigr)^n \hat C_n(k) \quad \text{if } n\geq2, \\
&g_n(k;z) = \Bigl( \frac z{1+2d\kappa D(1)} \Bigr)^n \hat\Pi_n(k) \quad \text{if } n\geq2.
\end{align}
\stopalpheqno

Assumptions S and D of
\cite{vHS} clearly hold, because of our assumptions \eqref{condha} and \eqref{condhb}, and
of the appendix of \cite{vHS}. There remains to check Assumption G by
using Proposition \ref{propestim}. In words, the task is to prove suitable bounds for $\|\Pi_n\|_1$,
assuming bounds for $\|C_n\|_\infty$. Let $\delta_0 = -\nabla^2 \hat D(0)$. The constant $\mu$
in the following lemma is any real number, not necessarily the connective constant.

\begin{lemma} $d\geq5$, and $\kappa$ satisfies the condition of Theorem \ref{thmconnconst}.
\label{lemestim}
\begin{itemize}
\item[(i)] Assume that $\|C_m\|_\infty \leq K \beta \mu^m m^{-d/2}$ and $\|C_m\|_1 \leq K \mu^m$
for all $m<n$. Then if $\beta$ is small enough, $\|\Pi_n\|_1 \leq K' \beta \mu^n n^{-d/2}$.
\item[(ii)] Assume in addition that $\sum_x |x|^2 C_m(x) \leq K \delta_0 m \mu^m$ for all $m<n$. Then
if $1\leq\gamma\leq2$ and $\beta$ is small enough,  $\sum_x |x|^{2\gamma} |\Pi_n(x)| \leq K'
\delta_0 \beta \mu^n n^{-d/2+\gamma}$.
\end{itemize}
The constant $K'$ depends only on $d,\kappa$, and $K$.
\end{lemma}

\begin{proof}
For item (i) we use Proposition \ref{propestim} (i) and (ii).
\ba
\|\Pi_n\|_1 \leq &(1+2d\kappa) K \beta \mu^{n-1} (n-1)^{-d/2} \nn\\
&+ \sum_{N\geq2} (2N-1) 2^{N-1} (1+2d\kappa)^N \sum \prodtwo{j=1}{\text{odd}}^{2N-1} K \beta
\mu^{m_j} m_j^{-d/2} \prodtwo{j=2}{\text{even}}^{2N-2} K \mu^{m_j}.
\end{align}
The unlabeled sum is as in Proposition \ref{propestim}. Since $m_1$ is the largest term, we have $m_1 \geq
\frac n{2N-1}$. The unlabeled sum is bounded by $K^{2N-1} (2N-1)^{d/2} \beta^N \mu^n n^{-d/2}$
multiplying
\be
\prodtwo{j=3}{\text{odd}}^{2N-1} \Bigl\{ \sum_{m_j \geq 1} m_j^{-d/2} \sum_{m_{j-1}=1}^{m_j} 1
\Bigr\} \leq \Bigl( \sum_{m\geq1} m^{-d/2+1} \Bigr)^{N-1}.
\end{equation}
Lemma \ref{lemestim} (i) is then clear, since the sum over $N$ converges for $\beta$ small enough,
and contributes less than $\const \cdot \beta \mu^n n^{-d/2}$.

For item (ii) we use Proposition \ref{propestim} (iii) and (iv). We first observe that under
the assumptions of Lemma \ref{lemestim} we have, for $m<n$,
\be
\label{boundsecmom}
\bigl\| |x|^2 C_m \bigr\|_\infty \leq K'' \delta_0 \beta \mu^m m^{-d/2+1}.
\end{equation}
Indeed, because of \eqref{subadd} we can write \cite{vHS2}
\ba
|x|^2 C_m(x) &\leq 2\sum_y (|y|^2 + |x-y|^2) C_{m/2}(y) C_{m/2}(x-y) \nn\\
&\leq 4 \| C_{m/2} \|_\infty \sum_y |y|^2 C_{m/2}(y).
\end{align}
(This was assuming $m$ to be even; the case $m$ odd is very similar.) Inserting
\eqref{boundsecmom} in Proposition \ref{propestim} (iv), one gets Lemma \ref{lemestim} (ii)
with $\gamma=2$. Then one can use H\"older's inequality to get \cite{vHS2}
\be
\sum_x |x|^{2\gamma} |\Pi_n(x)| \leq \Bigl( \sum_x |\Pi_n(x)| \Bigr)^{1-\gamma/2} \Bigl(
\sum_x |x|^4 |\Pi_n(x)| \Bigr)^{\gamma/2}.
\end{equation}
Item (ii) for general $\gamma$ follows from (i), and (ii) with $\gamma=2$.
\end{proof}

Assumption G of \cite{vHS} can be proved with the help of Lemma \ref{lemestim}.\footnote{The assumption of
\cite{vHS} involves a bound for $\|\hat D^2 f_m\|_1$ instead of $\|f_m\|_1$. The former
easily implies a suitable bound for $\|f_{m+2}\|_1$.} The first two inequalities are
straighforward, as is the third one since $\partial_z g_n(k;z) = \frac{1+2d\kappa D(1)}z n
g_n(k;z)$. The last inequality is more involved and deals with the error of the Taylor
expansion of $g_n$ to second order. First we write
\be
\hat\Pi_n(k) - \hat\Pi_n(0) - \frac{|k|^2}{2d} \nabla^2 \hat\Pi_n(0) = \sum_x \bigl[ \cos(kx)
- 1 + \tfrac12 (kx)^2 \bigr] \Pi_n(x).
\end{equation}
We used the symmetries of $\Pi_n(x)$ to replace a term $\frac1{2d} |k|^2 |x|^2$ by $\frac12 (kx)^2$. Now
$|\cos\xi - 1 + \frac12 \xi^2| \leq \const \, \xi^{2+\varepsilon}$ for any
$0\leq\varepsilon\leq2$. This, and Lemma \ref{lemestim}, clearly implies the validity of Assumption G for
$|k|$ small; when $|k|$ is large the situation is clear.

Theorem \ref{thmdiffbeh} is now an immediate consequence of Theorem 1.1 (b) of \cite{vHS}.

Finally, the diffusion constant $\delta$ can be given an explicit expression, once the convergence
of the lace expansion is established. See Eq.\ (3.5) of \cite{BR}, or Theorem 1.1
(d) of \cite{vHS}. The expression is
\be
\delta = \frac{\mu^{-1} \delta_0 + \tau}{1+\sigma},
\label{diffconst}
\end{equation}
where $\tau$ and $\sigma$ are given by
\ba
\tau &= -\sum_{m\geq1} \frac{\nabla^2 \, \hat\Pi_m(0)}{\mu^m}, \nn\\
\sigma &= \sum_{m\geq2} (m-1) \frac{\pi_m}{\mu^m}. \nn
\end{align}

\vspace{3mm}
\noindent
{\it Acknowledgments:}
It is a pleasure to thank Michael Aizenman, David Brydges, and Yvan Velenik for discussions
and encouragements. I am also indebted to Remco van der Hofstad and Gordon Slade for pointing out,
and clarifying it to me, that the diffusive behavior of the walk follows from their results
and Proposition \ref{propestim} here.

\end{document}